\documentclass{article}

\def \Q {{\bf Q}}

\def \P {{\bf P}}

\usepackage[T1]{fontenc} % РёРЅРѕРіРґР° T1 РІРјРµСЃС‚Рѕ T2A
\usepackage[cp1251]{inputenc} 
\usepackage[russian]{babel}

\begin{document}

  \rm \Large

\begin{center}   
{ \LARGE \bf
Mixing automorphisms not  having the Kolmogorov spectral  property
\vspace{5mm}

Valery V. Ryzhikov}
\end{center} 

\vspace{5mm}
 \bf Abstract.  \it  Let  $T$ be a staircase rank-one  construction with parameters $r_j \sim j^d$, $0<d<0.2$, then  its spectrum does not have the group property, and the product $T\odot T$ has  simple  spectrum. \rm

\section{Introduction}
 %  \bf  The group property of the spectrum. 
  \it Long ago Kolmogorov conjectured (see, for example, {\rm \cite{Si}}) that the maximum spectral type of an ergodic automorphism always subordinates its convolution. This property is a natural continual analogue of the group  property of the spectrum of an ergodic automorphism with a discrete  spectrum and was proved by Sinai for a special class of automorphisms  satisfying condition A (see {\rm \cite{Si}}). However, the conjecture is not true  in general. In this section we construct ergodic automorphisms whose  maximum spectral types do not have the "group property" in the above  sense. \rm  (cited from \cite{KS}).

\vspace{5mm}
The topic of the absence of a group property of the spectrum was developed: V.I. Oseledets constructed  weakly mixing examples \cite{O}, A.M. Stepin established the typicality of this property \cite{S}. These transformations, by virtue of their definition, are not  mixing. Mixing examples occurred due to another problem, known as the Rokhlin problem on non-simple homogeneous spectrum (for the history and evolution of this problem, see \cite{An}, \cite{Da}).

In \cite{99} the author proved the existence of  staircase constructions $T$ such that the products $T\otimes T$ have homogeneous spectrum of multiplicity 2. (Automorphism of the Lebesgue space $(X,\mu)$ and the unitary operator induced by it in $L_2 (X,\mu)$ are denoted identically.) In this case, its spectral measure $\sigma_T$ and its convolution square $\sigma_T^{\ast 2}$ are mutually singular. Taking into account the fact that T. Adams \cite{Ad} proved the mixing property of a suitable class of staircase  constructions, in \cite{99} two problems in the class of mixing automorphisms were solved at once: the absence of a group property of the spectrum for a mixing automorphism such that its  tensor square had a homogeneous spectrum of multiplicity 2. However, due to the specifics of the proof, based on “infinitesimal” weak limits, explicit examples were not indicated (as well as in the work \cite{Ag}, where it was used the  method \cite{99} of  approximations of mixing systems by non-mixing ones).
 
Specific automorphisms with these properties appeared in \cite{20}.
The staircase rank one construction  is specified by a sequence of parameters $r_j$. In \cite{20} it is shown that the necessary examples are  all constructions for which $r_j \sim \log j$. Let's expand the class of examples with appropriate explanations. \rm

\vspace{3mm}
\bf Theorem 1. \it If $T$ is a staircase construction with parameters $r_j \sim j^d$, $0<d<0.2$, then its spectrum does not have the group property, and the product $T\odot T$ has  simple  spectrum.\rm

\vspace{3mm}
The existence of mixing automorphisms with homogeneous spectrum of multiplicity $n>2$ is established in \cite{T}. There is no constructive solution yet.

\section{Staircase  constructions}
\it Definition of staircase transformation. \rm 
Let $r_j\to\infty$, and for all $j$, starting from $j>j_0$,  we set 
$\bar s_j=(1,2,\dots, r_j-2,r_j-1,0).$ 

Let at step $j$  a system of disjoint half-intervals
$$E_j, TE_j, T^2E_j,\dots, T^{h_j-1}E_j,$$
be given, and on half-intervals $E_j, TE_j, \dots, T^{h_j-2}E_j$
the transformation $T$ is a parallel translation. Such a set of half-intervals is called a tower of stage $j$; their union is denoted by $X_j$ and is also called a tower.

Let us represent $E_j$ as a disjoint union of $r_j$ half-intervals
$$E_j^1,E_j^2E_j^3,\dots E_j^{r_j}$$ of the same length.
For each $i=1,2,\dots, r_j$ we define the column $X_{i,j}$ as the union of intervals
$$E_j^i, TE_j^i ,T^2 E_j^i,\dots, T^{h_j-1}E_j^i.$$
To each column $X_{i,j}$ we add $i$ of disjoint half-intervals (spacers) of the same measure as $E_j^i$, obtaining a set
$$E_j^i, TE_j^i, T^2 E_j^i,\dots, T^{h_j-1}E_j^i, T^{h_j}E_j^i, T^{h_j+1}E_j^i , \dots, T^{h_j+s_j(i)-1}E_j^i$$
(all these sets do not intersect).
Denoting $E_{j+1}= E^1_j$, for $i<r_j$ we set
$$T^{h_j+s_j(i)}E_j^i = E_j^{i+1}.$$
 The set of superstructured columns is from now on considered as a tower of stage $j+1$, consisting of half-intervals
$$E_{j+1}, TE_{j+1}, T^2 E_{j+1},\dots, T^{h_{j+1}-1}E_{j+1},$$
where
 $$ h_{j+1} = \sum_{i=1}^{r_j}(h_j+s_j(i)).$$

The partial definition of the transformation $T$ at step $j$ is preserved in all subsequent steps. On the space $X=\cup_j X_j$, an invertible transformation $T:X\to X$ is thereby defined, preserving the standard Lebesgue measure on $X$.

This transformation $T$ is ergodic, has  simple spectrum. It is known that the indicators of $E_j$ are cyclic vectors for the operator $T$.

\section{Proof of Theorem 1}
We fix the indicator $f$ of some floor.
Let us prove that for the chosen staircase transformation  $T$ for all $r>0$ the vectors $T^rf\otimes f + f\otimes T^rf$
 belong to the cyclic space $C_{f\otimes f}$ of the operator $T\otimes T$. In fact, for this it is enough to show that  for the $L_2$-distance
$\rho$ we have  
$$\rho(T^rf\otimes f + f\otimes T^rf\,,\, C_{f\otimes f})\to 0.$$
Then the symmetric power $T\odot T$ has a simple spectrum, which implies the absence of the group property for $\sigma_T$.

Let us denote $$ Q_r = \frac 1 r\sum_{i=0}^{r-1} T^{-i}, \ \ \ \ \Q_r=Q_r\otimes Q_r.$$
We have
$$T^rf\otimes f + f\otimes T^rf
=r^2Q_rf\otimes Q_rf \ \ \ + \ \ \ (r-2)^2
 Q_{r-2}Tf\otimes Q_{r-2}Tf -$$
$$-(r-1)^2
 Q_{r-1}f\otimes Q_{r-1}f \
-\ (r-1)^2
 Q_{r-1}Tf\otimes Q_{r-1}Tf.\eqno (1)$$

Assuming the sequence $r_j\to\infty$ is monotonic, we set
$$ J_r=\{j:\, r_j=r+1, \ \ j<j_r-r\},$$
where $j_r=\max\{j:\, r_j=r+1,\}$.

We will find $D>0$ such that for $|J_r|> r^D$ we have
$$\left\|\Q_r(f\otimes f)- \P_r(f\otimes f)\right\|^2 <
\, \frac 2 {|J_r|}=o\left(\frac 1{r^4}\right), \eqno (2)$$
Where
$$ \P_r(f\otimes f) = \frac 1 {|J_r|}\sum_{j\in J_r} T^{h_j}f\otimes T^{h_j}f.$$
If $r\to\infty $ is satisfied
the norm of the difference $\P_r(f\otimes f)-\Q_r(f\otimes f)$ is $o\left(\frac 1{r^2}\right),$ then taking into account
equality $(1)$ we obtain that the vectors $T^rf\otimes {\bf 1} + {\bf 1}\otimes T^r$ "almost lie" in the cyclic space of the vector $f\otimes f$, which is what we need.

Consider the equality
$$ {|J_r|^2}\|\Q_r(f\otimes f)-\P_r(f\otimes f)\|^2 =\ {|J_r|^2}(\Q_r(f\otimes f) ,\Q_r(f\otimes f))- $$
$$ - 2{|J_r|}\left(\Q_r(f\otimes f)\,,\,\sum_{j\in J_r} T^{h_j}f\otimes T^{h_j}f\right) \ +\
\left(\sum_{j\in J_r} T^{h_j}f\otimes T^{h_j}f
\,,\,\sum_{k\in J_r} T^{h_k}f\otimes T^{h_k}f\right). \eqno (3)$$

If  $\left(\Q_r(f\otimes f)\,,\, T^{h_j}f\otimes T^{h_j}f\right)$ are sufficiently close to  $(\Q_r(f\otimes f),\Q_r(f\otimes f))$ and are close as well to 
$\left( T^{h_k}f\otimes T^{h_k}f\,,\,T^{h_m}f\otimes T^{h_m}f\right)$ for $k\neq m\in J_r$,
then the inequality $(2)$ we need will be satisfied. Now we say what means "close enough".

\vspace{3mm}
\bf Lemma. \it Let $j\in J_r$.
Then $$|(T^{h_j}f, Q_{r}f) - (Q_{r}f,Q_{r}f)|< \frac {2r}{h_j}+ 2r^{-r}= o(r^{-1000}). $$\rm

\vspace{3mm}
For $j,\in J_r$ the lemma implies
 $$|\left(\Q_r(f\otimes f)\,,\, T^{h_j}f\otimes T^{h_j}f\right)\, -\, (\Q_r(f\otimes f) ,\Q_r(f\otimes f))|=o(r^{-100}).$$
For some constant $C$ for $j, k=j+p\in J_r$, $p\geq 0$, we also have
$$|\left(T^{h_j}f\otimes T^{h_j}f\,,\, T^{h_{j+p}}f\otimes T^{h_{j+p}}f\right)\, -\, (\Q_r(f\otimes f),\Q_r(f\otimes f))|< \frac {Ch_j}{h_{j+p}}\leq Cr^{-p}. $$

From these inequalities $(2)$ follows, 

The main contribution (less than $\frac 2 {|J_r|}$) to the value of the difference, which is the right side of the equality $(3)$ ,  give 
$$\left(T^{h_j}f\otimes T^{h_j}f
\,,\,T^{h_{j+p}}f\otimes T^{h_{j+p}}f\right)
$$
as $p=0$. From this we can get $(2)$.
 Thus, for $d<\frac 1 5$ and $r_j\sim j^{d}$ we have $|J_r|\sim r^D\gg r^4$, where $D=\frac{1-d} d$, so
$$ \frac 1 {|J_r|}=o\left(\frac 1{r^4}\right).$$ 
This leads, as we explained, to a simple spectrum of the product $T\odot T$ for the corresponding  construction $T$.

In the case of infinite automorphisms (of a space with a sigma-finite measure, the corresponding “mixing”  automorphisms  ($T^n\to_w 0$)   without the Kolmogorov  property of spectrum are presented  among explicit Sidon’s constructions \cite{24}.

\end{document}